\input amstex
\TagsOnRight
\magnification\magstep1
\font\xivtn=cmr10 scaled \magstep2
\documentstyle{amsppt}
\parindent=12pt
\rightskip=0pt plus 2pt minus 2pt
\hsize16.5 true cm
\vsize23.0 true cm
\voffset=0 true cm
\hoffset=0 true cm
\parskip=2pt plus 3pt
\def\ord{\mathop{\fam0ord}\nolimits}
\let\epsilon\varepsilon
\let\phi\varphi

\noindent
{\eightrm Mat. Zapiski, {\eightbf 2}, Diophantine approximations,
Proceedings of papers dedicated to the memory of Prof. N.I.
Feldman, ed. Yu.V. Nesterenko, Centre for applied research
under Mech.-Math. Faculty of MSU, Moscow (1996), 23--42.
In russian: \hfill
http://www.math.jussieu.fr/~nesteren/
\hskip 3.4 true cm
\break
In english: \hfill
http://www.math.jussieu.fr/~miw/articles/ps/Nesterenko.ps
}

\bigskip
\bigskip
\bigskip

\centerline{{\xivtn On the approximation of the values of
exponential}}
\centerline{{\xivtn function and logarithm by algebraic numbers}}
\bigskip
\centerline{{\smc Yu.~Nesterenko (Moscow), M.~Waldschmidt (Paris)}}
\medskip
\line{\hfill{\sl Dedicated to the memory of Professor N.I.~Feldman}}
\bigskip
\noindent{\bf\S1. Introduction}

According to the theorem of Lindemann for any non-zero number
$\theta$ both of numbers $\theta$ and $e^\theta$ can not be algebraic.
For any algebraic numbers $\alpha$ and $\beta$ the expression
$|e^\theta-\alpha|+|\theta-\beta|$ does not vanish. How
small can this expression be? The answer should obviously depend
on the three following parameters: the heights
$h(\alpha)$, $h(\beta)$ of the algebraic numbers, the degree of the
number
field $\Bbb{D}=\Bbb{Q}(\alpha,\beta)$. Here, we denote by $h(\alpha)$
the
{\sl absolute logarithmic Weil height\/} of $\alpha$: when the
minimal polynomial
of $\alpha$ is $P(x)=a_0x^d+\cdots+a_{d}\in\Bbb{Z}[x]$, and its
complex
conjugates are $\alpha_1,\hdots,\alpha_d$:
$$
a_0x^d+\cdots+a_d=a_0(x-\alpha_1)\cdots(x-\alpha_d),
$$
then the absolute logarithmic Weil height $h(\alpha)$ of $\alpha$ is
defined by
$$
h(\alpha)=
\frac1d\left(\log|a_0|+\sum_{i=1}^d\log\max\{1,|\alpha_i|\}\right)
$$
while
$$
L(\alpha)=L(P)=\sum_{i=0}^d|a_i|
$$
is the length of the number $\alpha$ and of the polynomial $P$.
It is possible to prove that
$$
h(\alpha)\le d^{-1}\cdot\log L(\alpha)
\tag1.1
$$
(see [F~1982], Lemma~8.2).

In the next Theorem it will be more convenient to have a
parameter $E$, which will be choosen separately in each special
situation.

\noindent{\bf Theorem 1 (Main Theorem).} {\sl Let $\theta\in\Bbb{C}$,
$\theta\ne0$, and $\alpha,\beta$  be algebraic numbers; define
$\Bbb{K}=\Bbb{Q}(\alpha,\beta)$ and $D=[\Bbb{K}:\Bbb{Q}]$. Let $A$, $B$
and $E$ be positive real numbers with $E\ge e$ satisfying
$$
\log A\ge\max\bigl(h(\alpha),D^{-1}\bigr),\quad\log B\ge h(\beta).
$$
Then
$$
\gather
|e^{\theta}-\alpha|+|\theta-\beta|
\ge\exp\Bigl(-211D\bigl(\log B+\log\log A+4\log
D+2\log(E|\theta|_+)+10\bigr)\\
\cdot\bigl(D\log A+2E|\theta|+6\log E\bigr)
\cdot\bigl(3.3D\log(D+2)+\log E\bigr)
\cdot(\log E)^{-2}\Bigr),
\endgather
$$
where $|\theta|_+=\max(1,|\theta|)$}.

\medskip

From the inequality of our Main Theorem we deduce transcendence measures
for several numbers:
$\pi$, $\log2$, $e$ and  more generally $\log\alpha$ and $e^{\beta}$ (for
algebraic numbers
$\alpha$ and $\beta$, $\alpha\ne1$, $\beta\ne0$). A transcendence measure
of a transcendental complex number $\theta$ is a lower bound for
$|P(\theta)|$, when $P\in\Bbb{Z}[x]$ is a non-zero polynomial, in terms of
the degree of $P$ and of the length of $P$. For deducing the estimates of
the measure of transcendence we need the following assertion, connecting
the measure of transcendence and the measure of approximation by algebraic
numbers.

\noindent{\bf Lemma 1.} {\sl Let $\theta\in\Bbb{C}$. Assume that for any
algebraic number $\xi$ with $\deg\xi=d$ and $L(\xi)=L$, the inequality
$$
|\theta-\xi|\ge e^{-d\phi(d,L)}
$$
holds, where $\phi(x,y)$ is an increasing function of all arguments.
Then for any non-zero polynomial $P\in\Bbb{Z}[x]$ with $\deg P=N$ and
$L(P)=M$, we have}
$$
|P(\theta)|\ge e^{-d\phi(N,2^N M)}\cdot\bigl(4L\sqrt{N}\bigr)^{-N}.
$$

\noindent{\bf Proof.} See, for example, [F~1982], Lemma~3.7.

\noindent{\bf Theorem 2.} 1) {\sl Let $\xi$ be a real algebraic number,
$d=\deg\xi$, $L(\xi)\le L$, $L\ge3$. Then}
$$
|\pi-\xi|\ge\exp\bigl\{-1.2\cdot10^6d\cdot
(\log L+d\log d)\cdot(1+\log d)\bigr\}.
$$
2) {\sl If $P\in\Bbb{Z}[x]$, $P\ne0$, $\deg P\le d$, $L(P)\le L$, and
$L\ge3$, then}
$$
\bigl|P(\pi)\bigr|\ge\exp\bigl\{-2\cdot10^6d\cdot
(\log L+d\log d)\cdot(1+\log d)\bigr\}.
$$

For the proof of the first assertion we choose $\theta=\pi i$,
$\alpha=-1$,
$\beta=i\xi$, $E=e^2$, $\log A=D^{-1}$, $\log B=h(\xi)=h(\beta)$ and
note that
$D\le2d$. Since
$$
\gather
6.6d\log(2d+2)+\log E<11.2d(1+\log d),\\
d\bigl(h(\xi)+3\log(2d)+2\log\pi+14\bigr)
\le17(\log L+d\log d),\\
1+2E|\theta|+6\log E\le59.5
\endgather
$$
we derive the assertion 1).

The second assertion follows from the first one and Lemma~1.

By the same way can be proved

\noindent{\bf Theorem 3.} 1) {\sl Let $\xi$ be a real algebraic number
with $d=\deg\xi$, $L(\xi)\le L$ and $L\ge3$. Then}
$$
|\log2-\xi|\ge\exp\bigl\{-151000\cdot d^2\cdot(\log L+d\log d)
\cdot(1+\log d)^{-1}\bigr\}.
$$
2) {\sl If $P\in\Bbb{Z}[x]$, $P\ne0$, $\deg P\le d$, $L(P)\le L$, and
$L\ge3$, then}
$$
\bigl|P(\log2)\bigr|\ge\exp\bigl\{-2.6\cdot10^5d^2\cdot
(\log L+d\log d)\cdot(1+\log d)^{-1}\bigr\}.
$$

For the proof of the first assertion we choose $\theta=\log2$,
$\alpha=2$,
$\beta=\xi$, $E=eD$, $A=e$, $\log B=h(\beta)$. In this case $D=d$. We
deduce
$$
\gather
d\bigl(h(\xi)+4\log d+12\bigr)
\le13(\log L+d\log d),\\
3.3d\log(d+2)+\log(ed)<5d(1+\log d),\\
d+2E|\theta|+6\log E\le 11d.
\endgather
$$
Therefore the first inequality of the Theorem~3 holds. The second
one follows from the first and Lemma~1.

\noindent{\bf Theorem 4.} 1) {\sl Let $\xi$ be a real algebraic number
with $d=\deg\xi$, $L(\xi)\le L$, $L\ge3$. Then}
$$
|e-\xi|\ge\exp\bigl\{-76000\cdot d^2\cdot(\log L+d)\bigr\}.
$$
2) {\sl If $P\in\Bbb{Z}[x]$, $P\ne0$, $\deg P\le d$, $L(P)\le L$, and
$L\ge3$, then}
$$
\bigl|P(e)\bigr|\ge\exp\bigl\{-1.3\cdot10^5\cdot d^2\cdot(\log
L+d)\bigr\}.
$$

For the proof of this theorem we take $\theta=1$, $\alpha=\xi$,
$\beta=1$,
$\log A=1+d^{-1}\log L$, $B=1$, $E=ed\log A$, $D=d$. The desired
estimates follow from the inequalities
$$
\gather
3\log\log A+6\log d+12
\le9(1+\log D+\log\log A)=9\log E,\\
3.3d\log(d+2)+\log E\le\tfrac{10}{3}d\log E,\\
d\log A+2E|\theta|+6\log E\le12(d+\log L).
\endgather
$$

Taking $\theta=\beta$ or $\theta=\log\alpha$ for any determination of  the
logarithm of $\alpha$ we can prove a lower bound for $|e^\beta-\alpha|$ and
$|\log\alpha-\beta|$.

\noindent{\bf Theorem 5.} {\sl Let  $\alpha$ and $\beta$ be algebraic
numbers; define $\Bbb{K}=\Bbb{Q}(\alpha,\beta)$ and $D=[\Bbb{K}:\Bbb{Q}]$.
Let $A$ and $E$ be positive real numbers satisfying $E\ge e$ and
$$
\log A\ge\max\bigl(h(\alpha)\, ,\, D^{-1}\log E\,,\, D^{-1}|\beta| E\bigr).
$$
1) If $\beta\not=0$, then
$$
\gather
|e^{\beta}-\alpha|
\ge\exp\Bigl(-105500\cdot D^2\log A\cdot\bigl(h(\beta)+\log_+\log A+\log
D+\log E\bigr)\\
\cdot\bigl(D\log D+\log E\bigr)
\cdot(\log E)^{-2}\Bigr),
\endgather
$$
where $\log_+x=\log\max(1,x)$.}

\par\noindent
{\sl
2) If $\alpha\not=0$, and if $\log\alpha$ is any non-zero determination of
the logarithm of $\alpha$, then
$$
\gather
|\beta-\log\alpha|
\ge\exp\Bigl(-105500\cdot D^2\log A\cdot\bigl(h(\beta)+\log_+\log A+\log
D+\log E\bigr)\\
\cdot\bigl(D\log D+\log E\bigr)
\cdot(\log E)^{-2}\Bigr).
\endgather
$$
}

\medskip\noindent
1) For the proof of the first assertion we choose $\theta=\beta$ and we
use the estimates
$$
h(\beta)+\log\log A+4\log D+2\log(E|\beta|_+)+10
\le 12\bigl(h(\beta)+\log_+\log A+\log D+\log E\bigr),
$$
$$
D\log A+2E|\beta|+6\log E\le 9D\log A
$$
and
$$
9\cdot 12\bigl(3.3D\log(D+2)+\log E\bigr)\le 500(D\log D+\log E).
$$
2) The proof of the second assertion is essentially the same, with the
choice $\theta=\log\alpha$, using the estimate
$|\theta|\le|\beta|+|\beta-\theta|$.

\noindent{\bf Theorem 6.} 1) {\sl Let $\alpha$ be an algebraic number,
$\alpha\ne0,1$. Then there exists a constant $\gamma_1>0$, depending
only on $\alpha$ and the determination of the logarithm of $\alpha$
such that
if $P\in\Bbb{Z}[x]$, $P\not\equiv0$, $\deg P\le d$, $L(P)\le L$, then}
$$
\bigl|P(\log\alpha)|\ge\exp\bigl\{-\gamma_1d^2\cdot
(\log L+d\log d)\cdot(1+\log d)^{-1}\bigr\}.
$$
2) {\sl Let $\beta$ be an algebraic number, $\beta\ne0$. Then there
exists a constant $\gamma_2$, depending only on $\beta$, such that if
$P\in\Bbb{Z}[x]$, $P\not\equiv0$, $\deg P\le d$, $L(P)\le L$, then}
$$
\bigl|P(e^\beta)\bigr|\ge\exp\bigl\{-\gamma_2d^2\cdot(\log
L+d)\bigr\}.
$$

Theorem~6 follows from Theorem~5 with the help of Lemma~1.

There are plenty of results like our Theorems~2--6. The first
transcendence measure for the number $e$ goes back to Borel in 1899
[Bo~1899]. Early results on this subject, including works by Popken (1929)
and Mahler (1932), are quoted in [FS~1967]. We point out here that,
without explicit computation of the constants in the bounds, Theorem~2 was
proved for the first time by N.I.~Feldman, [F~1951,  1960]
and Theorem~6 by P.L.~Cijsouw [C~1974].  The main theorem of [D~1993]
provides a lower bound for
$|e^{\beta}-\alpha|$; the conclusion is that either the estimate of our
theorem 5  holds with the constant $105500$ replaced by $10^{11}$, or else
$$
|e^{\beta}-\alpha|
\ge e^{-10^{11}dD h(\beta)}
\quad\hbox{with}\quad d=[{\Bbb Q}(\beta):{\Bbb Q}].
$$
Further references are given in [W~1978] and [D~1993], as well as in
Feld'man's papers which are listed below.

\medskip

For the proof of the Main Theorem we use M.~Laurent's method
of interpolation determinants, which enables us to avoid the
construction of the auxiliary function and also to avoid the
extrapolation, to derive good constants in lower bounds. The organization
of this paper is as follows: in Section~2 we prove a variant of the zero
estimate of [LMN~1993]; Section~3 is devoted to analytic
estimates for Laurent's interpolation determinants. An important
tool in our proof is the use of binomial polynomials \S4. Next, in \S5,
we provide an arithmetic lower bound for non-zero algebraic
numbers (Liouville's inequality). The proof of the Main Theorem is
completed in Section~6.

\noindent{\bf\S2. Multiplicity estimate}
\par
\medskip
The proof of Hermite-Lindemann Theorem involves the complex
analytic functions $z$ and $e^{\beta z}$; for $P\in\Bbb{C}[X,Y]$,
the derivative $(d/dz)F$ of the function
$$
F(z)=P(z,e^{\beta z})
$$
is a polynomial in $z$ and $e^{\beta z}$, which we call $\delta P$:
$$
(d/dz)P(z,e^z)=\delta P(z,e^z).
$$
It is plain that $\delta$ is the derivative operator
$\frac{\partial}{\partial X}+\beta Y\frac{\partial}{\partial Y}$.
Hence we can define $\delta$ on $\Bbb{K}[X,Y]$ by
$$
\delta=\tfrac{\partial}{\partial X}+\beta Y\tfrac{\partial}{\partial
Y},
$$
when $\Bbb{K}$ is any field containing $\beta$. In this paper we work with
a field $\Bbb{K}$ of zero characteristic.

Here is our multiplicity estimate.

\noindent{\bf Lemma 2.} {\sl Let $\Bbb{K}$ be a field of zero
characteristic, $\beta$ a non-zero element of $\Bbb{K}$, and let $D_0$,
$D_1$, $S$ and $M$ be  positive integers  satisfying
$$
SM>(D_0+M)(D_1+1).
\tag2.1
$$
Let $(\xi_1,\eta_1),\hdots,(\xi_M,\eta_M)$ be elements in
$\Bbb{K}\times\Bbb{K}^*$ with $\xi_1,\hdots,\xi_M$
pairwise distinct. Then there is no non-zero polynomial
$P\in\Bbb{K}[X,Y]$,
of degree $\le D_0$ in $X$ and of degree $\le D_1$ in $Y$ which
satisfies}
$$
\delta^\sigma P(\xi_\mu,\eta_\mu)=0\quad
\text{{\sl for $1\le\mu\le M$ and $0\le\sigma<S$}.}
\tag2.2
$$

The proof is essentially the same as the proof of the zero
estimate in [LMN~1995]: we shall eliminate $Y$ using $D_1+1$
derivatives, and get a polynomial in $X$ which vanishes at $\xi_j$
with
multiplicity at least $S-D_1$.

\noindent{\bf Proof.} Let us suppose that a polynomial $P$ satisfies
all the
conditions of the lemma, equalities (2.2) and $P\ne0$. We assume, as
we may without loss of generality, that $Y$ does not divide the
polynomial $P$, and also that $P$ has degree $\ge1$ with respect to
$Y$.
Let us define the numbers $k_0=0<k_1<\hdots<k_n\le D_1$ by the conditions
$$
\gather
P(X,Y)=\sum_{i=0}^n Q_i(X)Y^{k_i},\\
Q_i(X)=b_iX^{m_i}+\cdots\in\Bbb{K}[X],\quad b_i\ne0,\quad
i=0,\hdots,n.
\endgather
$$

For $0\le\sigma\le n$, we consider the polynomials
$$
\delta^\sigma P(X,Y)=\sum_{i=0}^n Q_{\sigma i}(X)\cdot Y^{k_i},
\tag2.3
$$
where
$$
Q_{\sigma i}(X)=\sum_{j=0}^\sigma\binom\sigma{j}Q_i^{(\sigma-j)}(X)
(\beta k_i)^j=b_i(\beta k_i)^\sigma\cdot X^{m_i}+\cdots.
$$

It follows from this representation that the determinant
$$
\align
\Delta(X)
&=\det\bigl(Q_{\sigma i}(X)\bigr)_{0\le i,\sigma\le n}
=\det\bigl(b_i(\beta k_i)^\sigma\cdot X^{m_i}+\cdots\bigr)_{0\le
i,\sigma\le n}\\
&=b_0\hdots b_n\beta^{n(n+1)/2}\cdot B\cdot X^{m_0+\cdots+m_n}+\cdots,
\endalign
$$
where $B$ is a Vandermonde determinant constructed from the numbers
$k_0,\hdots,k_n$, hence $B\ne0$. Now from (2.3) we derive
$$
\Delta(X)=\sum_{\sigma=0}^n\Delta_\sigma(X,Y)\cdot\delta^\sigma
P(X,Y),
\quad\Delta_\sigma(X,Y)\in\Bbb{K}[X,Y],
$$
and for any $\tau\in\Bbb{Z}$, $0\le\tau<S-n$, with some
$c_{\tau,j,\sigma}\in\Bbb{K}$,
$$
\Delta^{(\tau)}(\xi_j)
=\sum_{\sigma=0}^{n+\tau}c_{\tau,j,\sigma}\cdot
\delta^\sigma P(\xi_j,\eta_j)=0,\quad j=1,\hdots,M.
$$

Since $n\le D_1$ and $\deg\Delta(X)=m_0+\cdots+m_n\le(n+1)D_0\le
D_0(D_1+1)$, we deduce
$$
(S-n)M\le\deg\Delta(X)\le D_0(D_1+1),
$$
and $SM\le D_0(D_1+1)+nM\le (D_0+M)(D_1+1)$. This contradicts to the
condition
(2.1) and completes the proof of lemma~2.

\medskip
\noindent{\bf\S3. Analytic upper bound}

We prove an upper bound for the absolute value of some
interpolation determinants; this estimate is a variant of some of
Laurent's results in [L~1989] and [L~1993].

\noindent{\bf Lemma 3.} {\sl Let $L$ be a positive integer, $E$, $M$, $S$,
and $\epsilon$ be positive real numbers with
$$
0<\epsilon<E^{-L}.
$$
For $1\le\lambda\le L$, let $b_{\lambda 1},\hdots,b_{\lambda L}$ be
complex
numbers, $\phi_\lambda(z)$ be a complex integral functions of one
variable;
further, for $1\le\mu\le L$, let $\zeta_\mu$ be a complex number and
$\sigma_\mu$ be a non-negative integer, $0\le\sigma_\mu\le S$. Assume
that for $1\le\lambda\le L$ and $1\le\mu\le L$ we have
$$
\log|b_{\lambda\mu}|\le M,\quad
\log\max_{z\le
E}\bigl|\phi_\lambda^{(\sigma_\mu)}(z\zeta_\mu)\bigr|\le M.
$$
Then the logarithm of the absolute value of the determinant
$$
\Cal{D}=\det\bigl\|\phi_\lambda^{(\sigma_\mu)}(\zeta_\mu)
+\epsilon b_{\lambda\mu}\bigr\|_{1\le\lambda,\mu\le L}
$$
is bounded by}
$$
L^{-1}\cdot\log|\Cal{D}|\le-\frac{L}{2}\cdot\log E+M+S\log E+\log(2LE).
$$

\noindent{\bf Proof.} Let us define
$$
a_{\lambda\mu}(z)=\phi_\lambda^{(\sigma_\mu)}(z\zeta_\mu)\quad\text{and}
\quad\Cal{D}(z)=\det\bigl\|a_{\lambda\mu}(z)+\epsilon
b_{\lambda\mu}\bigr\|
_{1\le\lambda,\mu\le L}.
$$
Then
$$
\Cal{D}(z)=\sum_{I\subset\{1,\hdots,L\}}
\epsilon^{L-|I|}\cdot\Cal{D}_I(z),
\tag3.1
$$
where
$$
\Cal{D}_I(z)=\det\bigl\|c_{\lambda\mu}(z)\bigr\|\quad\text{and}\quad
c_{\lambda\mu}(z)=\left\{\alignedat2
&a_{\lambda\mu}(z),\quad&&\text{if $\lambda\in I$},\\
&b_{\lambda\mu},&&\text{if $\lambda\notin I$}.
\endalignedat\right.
$$
We claim that the function of one variable $\Cal{D}_I(z)$ has a zero
at the
origin of multiplicity
$$
\ge\frac{|I|\cdot(|I|-1)}{2}-\sigma_1-\hdots-\sigma_L.
$$
The determinant $\Cal{D}_I(z)$ is a linear combination with constant
coefficients of determinants
$\Cal{D}_{I,J}(z)=\det\|a_{\lambda\mu}(z)\|_{\lambda\in I,\mu\in J}$,
where
$J$ runs all subsets of $\{1,\hdots,L\}$ with condition $|J|=|I|$.
For the proof of our claim it is sufficient to prove the inequality
$$
\ord\Cal{D}_{I,J}(z)\ge\frac{|I|\cdot(|I|-1)}{2}-\sigma_1-\hdots-
\sigma_L.
$$
By multilinearity we reduce the proof of this last inequality to
the special case $\phi_\lambda(z)=z^{n_\lambda}$ for some
$n_\lambda\in\Bbb{N}$, $\lambda\in I$. In this special case

$$
\align
\Cal{D}_{I,J}(z)&=\det\left(\binom{n_\lambda}{\sigma_\mu}\cdot\sigma_\mu!
\cdot(z\zeta_\mu)^{n_\lambda-\sigma_\mu}\right)_{\lambda\in I,\mu\in
J}\\
&=z^{\sum_{\lambda\in I}n_\lambda-\sum_{\mu\in J} \sigma_\mu}\cdot
\det\left(\binom{n_\lambda}{\sigma_\mu}\cdot\sigma_\mu!
\cdot\zeta_\mu^{n_\lambda-\sigma_\mu}\right)_{\lambda\in I,\mu\in J},
\endalign
$$
where the binomial coefficient $\binom{n_\lambda}{\sigma_\mu}$ means
0 if
$\sigma_\mu>n_\lambda$. If the right hand side is not identically
zero, then
the numbers $n_\lambda$, $\lambda\in I$, are pairwise distinct, and
then
the right hand side has a zero at the origin of multiplicity
$$
\ge\frac{|I|\cdot(|I|-1)}{2}-\sum_{\mu\in J}\sigma_\mu
\ge\frac{|I|\cdot(|I|-1)}{2}-\sigma_1-\hdots-\sigma_L.
$$
Our claim on the order of vanishing of $\Cal{D}_I(z)$ at the origin
easily
follows.

By means of the Schwarz lemma we conclude
$$
\align
\log|\Cal{D}_I(1)|
&\le-\biggl(\frac{|I|\cdot(|I|-1)}{2}
-\sigma_1-\hdots-\sigma_L\biggr)\log E+\log\max_{z\le E}
\bigl|\Cal{D}_I(z)\bigr|\\
&\le-\biggl(\frac{|I|\cdot(|I|-1)}{2}
-\sigma_1-\hdots-\sigma_L\biggr)\log E+L\log L+M|I|+M(L-|I|)\\
&\le-\frac{|I|^2}{2}\cdot\log E+\frac{|I|}{2}\log E
+ML+L\log L+SL\cdot\log E.
\endalign
$$
We derive now from (3.1)
$$
\align
\log|\Cal{D}|=\log|\Cal{D}(1)|
&\le L\log2-L^2\log E+ML+L\log L+SL\cdot\log E\\
&+\max_{I\subset\{1,\hdots,L\}}
\biggl(-\frac{|I|^2}{2}\cdot\log E+\bigl(L+\tfrac12\bigr)(\log E)
\cdot|I|\biggr).
\endalign
$$
The polynomial
$$
-\frac{\log E}{2}\cdot t^2+\bigl(L+\tfrac12\bigr)(\log E)\cdot t
$$
is an increasing function in the interval $1\le t\le L$. Then we have
$$
\align
L^{-1}\cdot\log|\Cal{D}|
&\le-\frac{L}{2}\log E+\bigl(L+\tfrac12\bigr)\log E
+\log(2L)-L\log E+M+S\log E\\
&<-\frac{L}{2}\cdot\log E+M+S\log E+\log(2LE).
\endalign
$$
This completes the proof of Lemma~3.

\medskip
\noindent{\bf\S4. Binomial polynomials}

\noindent
When $N$, $H$ be a non-negative integers, and $z$ a complex number,
let us define $\Delta(z,0,H)=1$, and
$$
\Delta(z,N,H)=
\biggl(\frac{z(z+1)\cdots(z+H-1)}{H!}\biggr)^q
\cdot\frac{z(z+1)\cdots(z+r-1)}{r!},
\tag4.1
$$
where
$$
N=qH+r,\quad 1\le r\le H.
$$
For $u$ a non-negative integer, we write $\Delta^{(u)}(z,N,H)$ for the
derivative $(d/dz)^u\Delta(z,N,H)$.

The first idea of eliminating the factorials from the
derivatives of auxiliary functions with the help of such polynomials
was introduced (in the case $H=N$) by Feldman in [F~1960,~a,b] for the
improvement of estimates of the measure of transcendence of $\pi$ and
logarithms of algebraic numbers. Later ([F~1968]) this was one of
his key tools in order to achieve a best possible dependence of
the estimate in terms of the heights of the coefficients $\beta_i$ in
lower bounds for linear combinations
$\beta_0+\beta_1\log\alpha_1+\cdots
+\beta_n\log\alpha_n$; in
turn, such an estimate has dramatic consequences, especially the
first effective improvement to Liouville's inequality. The
introduction of polynomials of this kind in the case $r=H$ and $H<N$
into the transcendence theory is due to A.~Baker [Ba~1972], who improved
in this way the dependence of lower bounds for linear forms in
logarithms in terms of the heights of the $\alpha_i$. The polynomials
(4.1) were introduced in [M~1994] where a more general
assertion than the next Lemma~4 is proved.

For each positive integer $k$ and real number $a$, we denote by $\nu(k)$
the least common multiple of $1,2,\hdots,k$ and by $[a]$~ the integer part
of $a$.

\noindent{\bf Lemma 4.} {\sl Let $N\ge1$, $H\ge1$, $\sigma\ge0$ and $x$
be integers. Define $d_\sigma=\nu(H)^\sigma$. Then
$$
d_\sigma\cdot\Delta^{(u)}(x,N,H)\in\Bbb{Z},\quad0\le u\le\sigma,
$$
and}
$$
\gather
\log d_\sigma<\tfrac{107}{103}\cdot\sigma H,\tag4.2\\
\sum_{u=0}^\sigma \binom\sigma{u}\cdot\bigl|\Delta^{(u)}(x,N,H)\bigr|
<\sigma^\sigma\cdot e^{N+H}\biggl(1+\frac{|x|}{H}\biggr)^N.\tag4.3
\endgather
$$

\noindent{\bf Proof.} Let $p$ be a prime number and $b_1\le
b_2\le\hdots
\le b_N$ be integers. For any integer $k>0$ we denote $r_k$
the number
of $b_i$ which are multiple of $p^k$. Then
$$
\ord_p(b_1\cdots b_N)=r_1+r_2+\cdots.
$$
If we delete any $u$ numbers from $b_1,\hdots,b_N$ and if
$b_{j_1},\hdots,
b_{j_{N-u}}$ denote the remaining $N-u$ numbers, we derive
$$
\ord_p(b_{j_1},\hdots,b_{j_{N-u}})\ge\sum_{k\ge1}\max(r_k-u,0).
$$
We define now the numbers $b_j$ as the $N$ factors in the product
$$
\bigl(x(x+1)\cdots(x+H-1)\bigr)^q\cdot
x(x+1)\cdots(x+r-1).
$$
In this case
$$
r_k\ge q\biggl[\frac{H}{p^k}\biggr]+\biggl[\frac{r}{p^k}\biggr],\quad
k\ge1.
$$
Now from the identity
$$
\Delta^{(u)}(z,N,H)=u!\cdot\Delta(z,N,H)\cdot
\sum(z+j_1)^{-1}\cdots(z+j_u)^{-1},
\tag4.4
$$
where summation is taken over all sets $\{j_1,\hdots,j_u\}$ such that
the polynomial $(z+j_1)\cdots(z+j_u)$ divides
$\Delta(z,N,H)$,
we see that
$$
\gather
\ord_p\bigl(d_\sigma\cdot\Delta^{(u)}(x,N,H)\biggr)
\ge\sigma\cdot\biggl[\frac{\log H}{\log p}\biggr]
-\sum_{p^k\le H}\Biggl(q\biggl[\frac{H}{p^k}\biggr]+
\biggl[\frac{r}{p^k}\biggr]\Biggr)\\
+\sum_{k\ge1}\max\Biggl(q\biggl[\frac{H}{p^k}\biggr]+
\biggl[\frac{r}{p^k}\biggr]-u,0\Biggr)
\ge\sum_{p^k\le H}\max\Biggl(\sigma-u,\sigma-
q\biggl[\frac{H}{p^k}\biggr]-
\biggl[\frac{r}{p^k}\biggr]\Biggr)\ge0.
\endgather
$$

This proves the assertion
$d_\sigma\cdot\Delta^{(u)}(x,N,H)\in\Bbb{Z}$,
$0\le u\le\sigma$, $x\in\Bbb{Z}$.

The estimate (4.1) follows from the inequality
$\log\nu(k)\le\tfrac{107}{103}
\cdot k$ (see for instance [Y~1989] Lemma~2.3 p.~127).

By the identity (4.4) we see
$$
\align
\sum_{u=0}^\sigma \binom\sigma{u}\cdot|\Delta^{(u)}(x,N,H)|
&\le\sum_{u=0}^\sigma \binom\sigma{u}\cdot \binom{N}{u}\cdot u!\cdot
\bigl(|x|+H-1\bigr)^{N-u}(H!)^{-q}(r!)^{-1}\\
&\le\sigma^\sigma\cdot\sum_{u=0}^N \binom{N}{u}\bigl(|x|+H-1\bigr)^{N-u}
(H!)^{-q}(r!)^{-1}\\
&\le\sigma^\sigma\bigl(|x|+H\bigr)^N\biggl(\frac{H^H}{H!}\biggr)^q
\biggl(\frac{H^r}{r!}\biggr)H^{-N}\\
&\le\sigma^\sigma\cdot e^{N+H}\cdot\biggl(1+\frac{|x|}{H}\biggr)^N.
\endalign
$$
This completes the proof of Lemma~4.

\medskip
\noindent{\bf\S5. Liouville's inequality}

For the next result, we use the notion of {\sl length\/} of a
polynomial $f\in\Bbb{C}[X_1,\hdots,X_n]$ (see \S 1).

\noindent{\bf Lemma 5 (Liouville's inequality).} {\sl Let $\Bbbk$ be a
subfield of $\Bbb{C}$ which is a finite extension of $\Bbb{Q}$ of
degree
$D$. Further let $\alpha_1,\hdots,\alpha_n$ be elements in $\Bbbk$.
Furthermore let $f$ be a polynomial in
$\Bbbk[X_1,\hdots,X_n]$, with coefficients in $\Bbb{Z}$, of degree at
most
$N_i$ with respect to $X_i$, and which does not vanish at the point
$(\alpha_1,\hdots,\alpha_n)$. Then
$$
\log\bigl|f(\alpha_1,\hdots,\alpha_n)\bigr|
\ge-(D'-1)\cdot\log L(f)-D'\sum_{i=1}^n N_ih(\alpha_i),
$$
where}
$$
D'=\left\{\alignedat2
&D/2\quad&&\text{{\sl if $\Bbbk$ is not a real field}},\\
&D\quad&&\text{{\sl if $\Bbbk$ is a real field}}.
\endalignedat\right.
$$

\noindent{\bf Proof.} See [F~1982], Lemma~9.2.

\medskip
\noindent{\bf\S6. Proof of the Main Theorem}

Let us suppose that under the conditions of the Main Theorem the inequality
$$
|e^\theta-\alpha|+|\theta-\beta|<E^{-211DUVW},
\tag6.1
$$
holds, with
$$
\gather
U=\frac{3.3D\log(D+2)+\log E}{\log E},\quad
V=\frac{2E|\theta|+D\log A+6\log E}{\log E},\\
W=\frac{\log B+\log\log A+4\log D+2\log(E|\theta|_+)+10}{\log E}.
\endgather
$$
Note that $U\ge1$, $V\ge6$, $W\ge2$.

a) {\sl Step one: Constuction of a non-zero determinant\/} $\Cal{D}$.

The proof of the Main theorem involves complex analytic
functions in one variable $\Delta(z,\tau,H)e^{\theta tz}$, for
non negative integers $\tau$, $H$ and $t$; the derivative of order
$\sigma$ of this function at the point $s\in\Bbb{Z}$, $s\ge0$, is
$$
\gamma_{\tau t}^{\sigma s}=\Bigl(\frac{d}{dz}\Bigr)^\sigma
\Bigl(\Delta(z,\tau,H)e^{\theta tz}\Bigr)\bigg|_{z=s}
=\sum_{k=0}^{\min(\tau,\sigma)}\frac{\sigma!}{(\sigma-k)!k!}\cdot
\Delta^{(k)}(s,\tau,H)\cdot(t\theta)^{\sigma-k}\cdot e^{\theta ts}.
\tag6.2
$$
We choose parameters $T$, $T_1$, $S$, $S_1$ and $H$:
$$
S=[10.5UV],\quad S_1=[12DW+0.5],
$$
$$
T=[20.2DVW],\quad
T_1=[4.2U+0.5],\quad H=[1.5W\log E]
$$
and restrict ourselves to the ranges
$$
0\le\tau\le T,\quad |t|\le T_1,\quad 0\le\sigma\le S,\quad |s|\le S_1.
$$
Replacing the numbers $e^\theta$ by $\alpha$ and $\theta$ by $\beta$
in (6.2),
we find an algebraic number
$$
\sum_{k=0}^{\min(\tau,\sigma)}\frac{\sigma!}{(\sigma-k)!k!}\cdot
\Delta^{(k)}(s,\tau,H)\cdot(t\beta)^{\sigma-k}\cdot\alpha^{ts},
$$
which will be a good approximation to $\gamma_{\tau t}^{\sigma s}$.
According to Lemma~4 the number
$$
a_{\tau t}^{\sigma s}=d_\sigma\cdot
\sum_{k=0}^{\min(\tau,\sigma)}\frac{\sigma!}{(\sigma-k)!k!}\cdot
\Delta^{(k)}(s,\tau,H)\cdot(t\beta)^{\sigma-k}\cdot\alpha^{ts},
\tag6.3
$$
will be an polynomial in $\alpha$, ${\alpha}^{-1}$, $\beta$ with
integer coefficients.

We also define $L=(T+1)(2T_1+1)$, which is the number of $(\tau,t)$.

\noindent{\bf Lemma 6.} {\sl There exists a set $\{(\sigma_\mu,s_\mu);
1\le\mu\le L\}$ of elements in $\Bbb{Z}\times\Bbb{Z}$ with
$0\le\sigma_\mu\le S$ and $0\le s_\mu\le S_1$ with the property that
the determinant of the $L\times L$ matrix
$$
\Cal{D}=\det\bigl\|a_{\tau}^{\sigma_\mu}{}_t^{s_\mu}\bigr\|,\quad
0\le\tau\le T,\quad |t|\le T_1,\quad 1\le\mu\le L,
$$
does not vanish}.

\noindent{\bf Proof.} Let $\Bbb{C}[X,Y,Y^{-1}]$ be the ring of
polynomials in $X$, $Y$, $Y^{-1}$
and let $\delta$ be the derivative operator on $\Bbb{C}[X,Y,Y^{-1}]$
defined by
$$
\delta=\tfrac{\partial}{\partial X}+\beta Y\tfrac{\partial}{\partial
Y}.
$$
Then
$$
\delta^\sigma\bigl(\Delta(X,\tau,H)Y^t\bigr)
=\sum_{k=0}^{\min(\tau,\sigma)}\frac{\sigma!}{(\sigma-k)!k!}\cdot
\Delta^{(k)}(X,\tau,H)\cdot(t\beta)^{\sigma-k}\cdot Y^t
$$
and
$$
a_{\tau t}^{\sigma s}=d_\sigma\delta^\sigma
\bigl(\Delta(X,\tau,H)Y^t\bigr)\Big|_{(X,Y)=(s,\alpha^s)}.
$$
Let us suppose that the rank of the matrix
$$
\|a_{\tau t}^{\sigma s}\|,\quad
0\le\tau\le T,\quad |t|\le T_1,\quad 0\le\sigma\le S,\quad |s|\le S_1,
$$
is less than $L$. Then there exist complex numbers $c_{\tau t}$,
$0\le\tau\le T$, $|t|\le T_1$, not all zero, such that the
polynomial
$$
R(X,Y)=\sum_{(\tau,t)}c_{\tau t}\Delta(X,\tau,H)Y^t \in \Bbb{C}[X,Y,Y^{-
1}]
$$
is not $0$ and satisfies
$$
\delta^\sigma R(X,Y)\Big|_{(X,Y)=(s,\alpha^s)}
=0,\quad 0\le\sigma\le S,\quad |s|\le S_1.
$$
But this contradicts  Lemma~2 with $P(X,Y)=Y^{T_1}R(X,Y)$, $D_0=T$,
$D_1=2T_1$,  $M=2S_1+1$, $\xi_s=s$, $\eta_s=\alpha^s$ and $S$ changed to
$S+1$: indeed, from the  inequalities
$$
2S_1+1\ge24DW,\quad 2T_1+1\le10.4U,\quad S+1\ge10.5UV,\quad
T\le20.2DVW,
$$
$V\ge6$ and
$$
\displaylines{\qquad
\frac{(T+2S_1+1)(2T_1+1)}{(S+1)(2S_1+1)}
=\frac{2T_1+1}{S+1}\cdot\biggl(\frac{T}{2S_1+1}+1\biggr)
<\frac{10.4}{10.5V}\cdot\biggl(\frac{20.2V}{24}+1\biggr)\hfill\cr\hfill
\le\frac{104}{105}\cdot\biggl(\frac{101}{120}+\frac{1}{6}\biggr)<1
\qquad\cr}
$$
we derive $P=0$. This completes the proof of Lemma~6.

b) {\sl upper bound for\/} $|\Cal{D}|$.

We plan to use Lemma~3 with $\lambda$ replaced by $(\tau,t)$, for the $L$
functions
$$
f_{\tau t}(z)=\Delta(z,\tau,H)\cdot e^{\theta tz},\quad
0\le\tau\le T,\quad |t|\le T_1,
$$
with the points $\zeta_\mu=s_\mu$, $1\le\mu\le L$, with
$\epsilon=E^{-211DUVW}$ and $b_{\tau t\mu}$, instead
$b_{\lambda\mu}$ in Lemma~3, defined by
$$
d_\sigma^{-1}\cdot a_\tau^{\sigma_\mu}{}_t^{s_\mu}
=\gamma_\tau^{\sigma_\mu}{}_t^{s_\mu}+\epsilon\cdot b_{\tau t\mu}.
$$
The estimates
$$
T+1\le20.2DVW+1\le\bigl(20.2+\tfrac{1}{12}\bigr)DVW,\quad
T_1+\tfrac{1}{2}\le5.2U,
$$
means that
$$
L=(T+1)(2T_1+1)<211DUVW
\tag6.4
$$
and $\epsilon<E^{-L}$.

{}From Lemma~4 and (6.2) we deduce
$$
\align
\max_{|z|\le E}\bigl|f_{\tau t}^{(\sigma)}(sz)\bigr|
&\le\max_{|z|\le
E}\sum_{k=0}^{\min(\tau,\sigma)}\frac{\sigma!}{(\sigma-
k)!k!}\bigl|\Delta^{(k)}(sz,\tau,H)\bigr|\cdot|t\theta|^{\sigma-k}
e^{|\theta zts|}\\
&\le S^S\cdot e^{H+T}\cdot\biggl(1+\frac{ES_1}{H}\biggr)^T\cdot
\bigl(|\theta|_+T_1\bigr)^S\cdot
e^{|\theta|ET_1S_1}\le e^M,
\endalign
$$
where
$$
M=S\log S+H+T+T\log\biggl(1+\frac{ES_1}{H}\biggr)
+S\log\bigl(E|\theta|_+T_1\bigr)+E|\theta|S_1T_1+S_1T_1.
$$
It follows from (6.1) that
$$
\gather
\max\bigl(|\beta|,|\theta|)\le|\theta|_+(1+(2S)^{-1}),\\
\max\bigl(|\alpha|,|\alpha|^{-1},|e^\theta|,|e^{-
\theta}|\bigr)\le e^{|\theta|}\cdot
\bigl(1+(2S_1T_1+2)^{-1}\bigr),
\endgather
$$
and for any integer $k$ and $\ell$ with $0\le k\le S$ and $|\ell|\le
S_1T_1$,
$$
\eqalign{
|\beta^k\alpha^\ell-\theta^ke^{\theta \ell}|
&\le|\beta|^k\cdot|\alpha^\ell-e^{\theta
\ell}|+|e^\theta|^\ell\cdot|\beta^k-\theta^k|\cr
&\le\epsilon\cdot|\theta|_+^k\cdot e^{(|\ell|+1)|\theta|}\cdot
\biggl(1+\frac{1}{2S}\biggr)^k\cdot
\biggl(1+\frac{1}{2S_1T_1+2}\biggr)^{|\ell|+1}\cdot\max(|\ell|,k)\cr
&\le\epsilon
e|\theta|_+^S\cdot e^{2S_1T_1|\theta|}\cdot\max(S,S_1T_1).
\cr}
$$
Now we use the inequalities
$$
eS\le e^S\le E^S,\quad
eS_1T_1\le e^{S_1T_1}
$$
and we write $\sigma$ and $s$ in place of $\sigma_\mu$ and $s_\mu$. From
(6.2), (6.3) and Lemma~4 we derive
$$
\align
\epsilon\cdot|b_{\tau t\mu}|
&=|d_\sigma^{-1}\cdot a_{\tau t}^{\sigma s}-\gamma_{\tau t}^{\sigma
s}|\\
&\le\sum_{k=0}^{\min(\tau,\sigma)}\frac{\sigma!}{(\sigma-k)!k!}\cdot
\bigl|\Delta^{(k)}(s,\tau,H)\bigr|\cdot |t|^{\sigma-k}\cdot
|\beta^{\sigma-k}\alpha^{ts}-\theta^{\sigma-k}e^{\theta ts}|\\
&\le S^S\cdot e^{H+T+1}\cdot\biggl(1+\frac{S_1}{H}\biggr)^T\cdot
T_1^S\cdot
|\theta|_+^S e^{2|\theta|T_1S_1}\max(S,S_1T_1)\cdot\epsilon
\le\epsilon\cdot e^M.
\endalign
$$
Since $\log d_\sigma\le\tfrac{107}{103}SH$ we deduce from Lemma~3
$$
\align
\frac{1}{L}\cdot\log|\Cal{D}|
&\le-\frac{L}{2}\cdot\log E+\tfrac{107}{103}SH+S\log S+H+T
+T\log\biggl(1+\frac{ES_1}{H}\biggr)\\
&+S\log\bigl(E|\theta|_+T_1\bigr)+|\theta|ES_1T_1+S_1T_1
+S\log E+\log(2LE).
\tag6.5
\endalign
$$

c) {\sl lower bound for\/} $|\Cal{D}|$.

Let us define the polynomials
$$
q_{\tau t}^{\sigma s}(X,Y)
=d_\sigma\cdot\sum_{k=0}^{\min(\tau,\sigma)}\frac{\sigma!}{(\sigma -
k)!k!}\cdot
\Delta^{(k)}(s,\tau,H)\cdot(tX)^{\sigma-k}\cdot Y^{ts}
\in\Bbb{Z}[X,Y,Y^{-1}],
$$
and
$$
R(X,Y)=\det\bigl\|q_\tau^{\sigma_\mu}{}_t^{s_\mu}(X,Y)\bigr\|,
\quad 0\le\tau\le T,\quad |t|\le T_1,\quad 1\le\mu\le L.
$$
It follows from (6.3) that
$a_{\tau t}^{\sigma s}=q_{\tau t}^{\sigma s}(\beta,\alpha)$
and $\Cal{D}=R(\beta,\alpha)$.  From  the inequalities
$$
\deg_X q_{\tau t}^{\sigma s}(X,Y)\le S,
\quad\deg_Y q_{\tau t}^{\sigma s}(X,Y)\le |t|S_1,
\quad\deg_{Y^{-1}} q_{\tau t}^{\sigma s}(X,Y)\le |t|S_1
$$
we derive
$$
\deg_X R(X,Y)\le LS,
$$
$$
\deg_Y R(X,Y)\le \tfrac{1}{4}LS_1(T_1+0.5),
\quad\deg_{Y^{-1}} R(X,Y)\le \tfrac{1}{4}LS_1(T_1+0.5).
$$
It follows from Lemma~4 that
$$
L\bigl(q_{\tau t}^{\sigma s}(X,Y)\bigr)
\le\exp\bigl(\tfrac{107}{103}SH\bigr)\cdot S^S\cdot e^{H+T}\cdot
\biggl(1+\frac{S_1}{H}\biggr)^T\cdot T_1^S
$$
and
$$
L\bigl(R(X,Y)\bigr)
\le L^L\cdot\Biggl(\exp\bigl(\tfrac{107}{103}SH\bigr)\cdot S^S\cdot
e^{H+T}\cdot
\biggl(1+\frac{S_1}{H}\biggr)^T\cdot T_1^S\Biggr)^L.
$$
Next from Lemma~5 we derive
$$
\align
\frac1L\cdot\log|\Cal{D}|
&\ge-(D-1)\Biggl(\log L+\tfrac{107}{103}SH+S\log S+H+T
+T\log\biggl(1+\frac{S_1}{H}\biggr)\\
&+S\log T_1\Biggr)-DS\log B-\tfrac12DS_1(T_1+0.5)\log A.
\tag6.6
\endalign
$$

d) {\sl End of the proof of Theorem 1}.

Let us compare the upper bound (6.5) for $|\Cal{D}|$ and the lower bound
(6.6). We derive
$$
\gather
\frac{L}{2}\log E\le\tfrac12S_1(T_1+0.5)
\Bigl(D\log A+2E|\theta|+2\Bigr)
+\Biggl(DT\log\biggl(1+\frac{S_1}{H}\biggr)+DT+T\log E\Biggr)\\
+DS\Bigl(\log B+\log S+\log\bigl(E|\theta|_+T_1\bigr)\Bigr)+DH
+DSH\tfrac{107}{103}+S\log E+\log(2E)+D\log L.
\endgather
$$
Using the definition of all parameters we deduce
$$
\tfrac12S_1(T_1+0.5)\Bigl(D\log A+2E|\theta|+2\Bigr)
\le0.5\cdot12.25\cdot5.2DUVW\log E=31.85DUVW\log E,
\tag6.7
$$
$$
\gather
\log\biggl(1+\frac{S_1}{H}\biggr)
\le\log\biggl(1+12.25\frac{D}{\log E}\biggr)
\le\log13.25+\log D\le2.6+\log D,
\endgather
$$
$$
\gather
DT\log\biggl(1+\frac{S_1}{H}\biggr)+DT+T\log E
\le20.2(D^2VW\log D+3.6D^2VW+DVW\log E)\\
\le20.2DVW\bigl(\log E+3.3D\log(D+2)\bigr)\le20.2DUVW\log E,
\tag6.8
\endgather
$$
$$
U\le1+3.3D\log(D+2)\le4.7D^{3/2},\quad
V\le9.8E|\theta|_+D\log A,
\tag6.9
$$
$$
\gather
\log\bigl(ST_1E|\theta|_+\bigr)
\le\log\bigl(50U^2VE|\theta|_+\bigr)
\le\log\log A+4\log D+2\log\bigl(E|\theta|_+\bigr)+10,\\
DS\Bigl(\log B+\log S+\log\bigl(E|\theta|_+T_1\bigr)\Bigr)
\le10.5DUV\Bigl(\log B+\log\log A+4\log D\\
+2\log\bigl(E|\theta|_+\bigr)+10\Bigr)
\le10.5DUVW\log E
\tag6.10
\endgather
$$
The estimates (6.4) and (6.9) mean that
$$
\align
\log L&\le\log(211DUVW)\le10+3.5\log D+\log\bigl(E|\theta|_+\bigr)
+\log\log A+\log W\\
&\le W\log E+\log W\le1.4W\log E\le 0.24UVW\log E.
\endalign
$$
With the help of this estimate and inequalities
$$
\align
DSH\le15.75DUVW\log E,\quad&
DH\le1.5DW\log E\le0.25DUVW\log E,\\
S\log E\le10.5UV\log E\le5.25DUVW\log E,\quad&
\log(2E)\le2\log E\le\tfrac16DUVW\log E,
\endalign
$$
we derive
$$
\displaylines{
\qquad DH+\tfrac{107}{103}DSH+S\log E+\log(2E)+D\log L\hfill\cr\hfill
\le\bigl(15.75\cdot\tfrac{107}{103}+5.25+0.49+\tfrac16\bigr)DUVW\log E
<22.28DUVW\log E.\quad(6.11)
\cr}
$$
Finally from (6.6)--(6.11) we deduce
$$
\eqalign{
\frac{L}{2}\log E&\le(31.85+20.2+10.5+22.28)DUVW\log E\cr
&= 84.83DUVW\log E<\frac{(T+1)(2T_1+1)}{2}\cdot\log E=\frac{L}{2}\log E.\cr}
$$
This contradiction means that (6.1) is wrong and completes the
proof of the Theorem 1.

\medskip
\noindent{\bf References}

\parindent=0pt

\medskip
[Ba~1972] A. Baker. A sharpening of the bounds for linear forms in
logarithms; Acta Arith., 21
(1972), 117--129.

\medskip
[Bo~1899] E.~Borel,
Sur la nature arithm\'etique du nombre $e$;
C.R.~Acad Sci.~Paris, 128 (1899), 596--599.

\medskip
[C~1974] P.L.~Cijsouw,
Transcendence measure of exponentials and logarithms of algebraic
numbers;
Compositio Math., 28 (1974), 163--178.

\medskip
[D~1993] G.~Diaz,
Une nouvelle minoration de $|\log\alpha-\beta|$,
$|\alpha-\exp\beta|$, $\alpha$ et $\beta$ alg\'ebriques;
Acta Arith., 64 (1993), 43--57;
compl\'ements: S\'eminaire d'Arithm\'etique Univ.~St.~Etienne 1990--91--92,
N$^\circ$4, 49--58.

\medskip
[F~1949] N.I.~Feldman,
The approximation of some transcendental numbers;
Dokl.~Akad.\ Nauk SSSR, 66 (1949), 565--567.

\medskip
[F~1951] N.I.~Feldman,
Approximation of certain transcendental numbers, I;
Izv.~Akad.~Na\-uk SSSR, Ser.~Mat., 15 (1951), 53--74;
engl.~transl.: Amer.~Math.~Soc.~Transl. (2) 59 (1966), 224--245.

\medskip
[F~1960a] N.I.~Feldman,
On the measure of transcendence of $\pi$;
Izv.~Akad.~Nauk SSSR, Ser.~Mat., 24 (1960), 357--368;
engl.~transl.: Amer.~Math.~Soc.~Transl. (2) 58 (1966), 110--124.

\medskip
[F~1960b] N.I.~Feldman,
Approximation of the logarithms of algebraic numbers by algebraic
numbers;
Izv.~Akad.~Nauk SSSR, Ser.~Mat., 24 (1960), 475--492;
engl.~transl.: Amer.~Math.~Soc.~Transl. (2) 58 (1966), 125--142.

\medskip
[F~1963a] N.I.~Feldman,
On the problem of the measure of transcendence of $e$;
Usp.~Mat.~Na\-uk 18:3 (111) (1963), 207--213.

\medskip
[F~1963b] N.I.~Feldman,
On the problem of the transcendence measure of $\pi$;
Usp.~Mat.~Na\-uk 18 (1963), 207--213.

\medskip
[F~1968] N.I.~Feldman,
Improved estimates for a linear form of the logarithms of algebraic
numbers;
Mat.~Sb., 77 (1968), 423--436;
engl.~transl.: Math.~USSR Sb., 77 (1968), 393--406.

\medskip
[F~1977] N.I.~Feldman,
Approximation of number $\pi$ by algebraic numbers from special
fields;
J.~Number Theory, 9 (1977), 48--60.

\medskip
[F~1982] N.I.~Feldman,
The seventh Hilbert's problem;
Moscow, Moscow State University, 1982, pp.~1--311.

\medskip
[FS~1967] N.I.~Feldman and A.B.~Shidlovskii,
The development and present state of the theory of transcendental
numbers;
Usp.~Mat.~Nauk, 22 (1967), 3--81;
engl.~transl.: Russian Math.~Surveys, 22 (1967), 1--79.

\medskip [L~1989]
M.~Laurent, Sur quelques r\'esultats r\'ecents de transcendance;
Journ\'ees arith\-m\'etiques Luminy, 1989, Ast\'erisque, 198--200 (1991),
209- -230.

\medskip
[L~1993] M.~Laurent,
Linear forms in two logarithms and interpolation determinants;
Acta Arith{.} {\bf 66} (1994), no{.} 2, 181--199.

\medskip
[LMN~1995] M.~Laurent, M.~Mignotte et Yu.~Nesterenko,
Formes lin\'eaires en deux logarithmes et d\'eterminants d'interpolation;
 J{.} Number Theory {\bf 55} (1995), no{.}
2, 285--321.

\medskip
[M 1993] E.M. Matveev. On the arithmetic properties of the values of
generalized binomial coefficients; Mat. Zam., 54 (1993), 76--81; Math.
Notes, 54 (1993), 1031--1036.

\medskip
[W~1978] M.~Waldschmidt,
Transcendence measures for numbers connected with the exponential
function;
J.~Austral.~Math.~Soc., 25 (1978), 445--465.

\medskip
[Y~1989] Yu Kunrui,
Linear forms in $p$-adic logarithms;
Acta Arith., 53 (1989), 107--186; (II), Compositio Math., 74 (1990),
15--113.

\vskip 2 true cm
\bigskip
\line{\hfill
\vtop{\hbox{Yu. V. NESTERENKO}
      \hbox{Moscow University}
      \hbox{Department of Mathematics and Mechanics}
      \hbox{119899 MOSCOW}
      \hbox{RUSSIA}
      \hbox{}
      \hbox{}
      \hbox{\tt NEST\@nw.math.msu.su}
      }\hfill
\vtop{\hbox{M. WALDSCHMIDT}
      \hbox{Universit\'e P. et M. Curie}
      \hbox{Institut de Math\'ematiques de Jussieu}
      \hbox{4, Place Jussieu, Case 247}
      \hbox{75252 PARIS Cedex 05}
      \hbox{FRANCE}
      \hbox{}
      \hbox{\tt miw\@math.jussieu.fr}
}\hfill}

\vfill
\noindent
{\eightrm
\hfill
http://www.math.jussieu.fr/~miw/articles/ps/Nesterenko.mai95.ps}
\bigskip
\end